\documentclass{jnmp}

\JNMPnumberwithin{equation}{section}

\usepackage{amsmath,amscd}

\newtheorem{theorem}{Theorem}
\newtheorem{proposition}{Proposition}

\begin{document}

\renewcommand{\evenhead}{H~Gargoubi}
\renewcommand{\oddhead}{Algebra ${\rm gl}(\lambda)$ Inside the Algebra
of Differential Operators on the Real Line}

\thispagestyle{empty}

\FirstPageHead{9}{3}{2002}{\pageref{Gargoubi-firstpage}--\pageref{Gargoubi-lastpage}}{Letter}

\copyrightnote{2002}{H~Gargoubi}

\Name{Algebra $\boldsymbol{{\mathrm{gl}}(\lambda)}$ Inside the Algebra \\
of Differential Operators on the Real Line}\label{Gargoubi-firstpage}

\Author{H~GARGOUBI}

\Address{I.P.E.I.M., route de Kairouan, 5019 Monastir, Tunisia\\
E-mail: hichem.gargoubi@ipeim.rnu.tn}

\Date{Received November 22, 2001; Revised February 26, 2002;
Accepted March 12, 2002}

\begin{abstract}
\noindent
The Lie algebra ${\rm gl}(\lambda)$ with $\lambda \in {\mathbb C}$,
introduced by B~L~Feigin, can be embedded into the Lie algebra of differential operators
on the real line (see~\cite{fe}). We give an explicit formula of the embedding of
${\rm gl}(\lambda)$
into the algebra ${\cal D}_{\lambda}$ of differential operators on the space of
tensor densities of degree $\lambda$ on ${\mathbb R}$.
Our main tool is the notion of projectively equivariant symbol
of a differential operator.
\end{abstract}

\section{Introduction}

The Lie algebra ${\rm gl}(\lambda)$ ($\lambda \in {\mathbb C}$) was introduced by
B~L~Feigin in~\cite {fe}
for calculation the cohomology of the Lie algebra of
differential operators on the real line.
The algebra ${\rm gl}(\lambda)$ is defined as the quotient of
the universal enveloping algebra $\mathrm{U}(\mathrm{sl}_2)$ of $\mathrm{sl}_2$
with respect to the ideal generated by the element
$\Delta - \lambda (\lambda - 1)$, where $\Delta$
is the Casimir element of $\mathrm{U}(\mathrm{sl}_2)$.
${\rm gl}(\lambda)$ is turned into a Lie algebra by the standard method of setting
$[a,b] =  ab  -ba$.

According to Feigin, ${\rm gl}(\lambda)$ can be considered as an analogue of ${\rm gl}(n)$
for $n = \lambda \in {{\mathbb N}}$; it is also called the algebra of
matrices of complex size, see also \cite{km,ls,sh,gl}.

We consider the space ${\cal D}_{\lambda}$ of all linear differential operators acting on
tensor densities of degree $\lambda$ on ${{\mathbb R}}$. One of the main results of \cite{fe} is
the construction of an embedding ${\rm gl}(\lambda) \to {\cal D}_{\lambda}$.

The purpose of this paper is to give an explicit formula of this embedding. We
also show that this embedding realizes  the isomorphism of Lie algebras
${\rm gl}(\lambda) \cong {\cal
D}^{\rm pol}_{\lambda}$ constructed in
\cite {bb1,bb2}, where ${\cal D}^{\rm pol}_{\lambda}
\subset {\cal D}_{\lambda}$ is the subalgebra of differential
operators with polynomial coefficients.

The main idea of this paper is to use the {\it projectively
equivariant symbol} of a differential operator, that is an
$\mathrm{sl}_2$-equivariant way to associate a polynomial function
on $T^*{\mathbb R}$ to a~differential operator. The notion of
projectively equivariant symbol was defined in \cite{cmz,lo} and
used in \cite{ga,go,go1} for study of modules of differential
operators.

\section{Basic definitions}

{\bf 2.1 The Lie algebra $\boldsymbol{{\rm gl}(\lambda)}$.}
Let $\mathrm{Vect}({\mathbb R})$ be the Lie algebra of smooth vector fields on~${\mathbb R}$
with complex coefficients: $X=X(x)\partial$, where $X(x)$ is a smooth complex
function of one real variable;
$X(x) \in C^{\infty}({\mathbb R},{\mathbb C})$, and
where $\partial=\frac{d}{dx}$.
Consider the Lie algebra $\mathrm{sl}_2\subset\mathrm{Vect}({\mathbb R})$
generated by the vector fields
\begin{equation}
\left\{\partial,x\partial,x^2\partial\right\}.
\label{gensl}
\end{equation}
Denote $ e_i := x^i\partial$, $i=0,1,2$,
the Casimir element
\[
\Delta
:=
e^2_1
-
\frac{1}{2} (e_0 e_2+ e_2 e_0)
\]
generates the center of $\mathrm{U}(\mathrm{sl}_2)$.
The quotient
\[
{\rm gl}(\lambda)
:=
\mathrm{U}(\mathrm{sl}_2) / (\Delta - \lambda (\lambda - 1)), \qquad \lambda \in {{\mathbb C}}
\]
is naturally a Lie algebra containing $\mathrm{sl}_2$.

{\bf 2.2 Modules of differential operators on $\boldsymbol{{\mathbb R}}$.}
Denote ${\cal D}$ the Lie algebra of linear differential
operators on ${{\mathbb R}}$ with complex coefficients:
\begin{equation}
A = a_n(x){\partial}^n+a_{n-1}(x){\partial}^{n-1}+\cdots+a_0(x),
\label{ope1}
\end{equation}
with $a_i(x) \in C^{\infty}({\mathbb R},{\mathbb C})$.

For $\lambda \in {\mathbb C}$,
$\mathrm{Vect}({\mathbb R})$ is embedded into the Lie algebra ${\cal D}$ by:
\begin{equation}
X \mapsto L_X^{\lambda}:=X(x)\partial+\lambda X^{\prime}(x) .
\label{lieder}
\end{equation}
Denote ${\cal D}_{\lambda}$
the $\mathrm{Vect}({\mathbb R})$-module structure
with respect to the adjoint action of $\mathrm{Vect}({\mathbb R})$ on~${\cal D}$.
The module ${\cal D}_{\lambda}$
has a natural filtration:
${\cal D}^0_{\lambda}\subset{\cal D}^1_{\lambda}\subset
\cdots\subset{\cal D}^n_{\lambda}\subset\cdots$,
where ${\cal D}^n_{\lambda}$ is the module of $n$-th order differential operators~(\ref{ope1}).

Geometrically speaking, differential operators are
acting on tensor densities, namely:
$A:{\cal F}_\lambda \to {\cal F}_\lambda$,
where ${\cal F}_\lambda$ is the space of
tensor densities of degree $\lambda$ on ${\mathbb R}$ (i.e., of sections of the line
bundle $(T^{\ast}{\mathbb R})^{\otimes\,\lambda}, \lambda\in {\mathbb C}$), that is:
$\phi =\phi (x)(dx)^{\lambda }$, where $\phi (x) \in C^{\infty}({\mathbb R},{\mathbb C})$.

It is evident that ${\cal F}_\lambda \cong C^{\infty}({\mathbb R},{\mathbb C})$
as linear spaces (but not as modules)
for any $\lambda$. We use this identification  throughout this paper.
The Lie algebra structures
of differential operators acting on the space of tensor densities
and on the space of functions are also identified  (see~\cite{ga}).

The $\mathrm{Vect}({\mathbb R})$-modules ${\cal D}_{\lambda}$
were considered by classics (see \cite{car,wil})
and, recently, stu\-died in a series of papers
\cite{do,go,ga,go1,lmt}.

{\bf 2.3 Principal symbol.}
Let  $\mathrm{Pol}(T^*{{\mathbb R}})$ be the space of functions on $T^*{{\mathbb R}}$
polynomial in the fibers. This space is usually considered as the space of
symbols associated to the space of differential operators on ${{\mathbb R}}$.

Recall that the {\it principal symbol} of a differential operator
is the linear map
$\sigma :{{\cal D} }\to \mathrm{Pol}(T^*{{\mathbb R}})$
defined by:
\[
\sigma (A) = a_n(x) \xi^n  ,
\]
where $A$ is a differential operator (\ref{ope1}) and
$\xi$ is the coordinate on the fiber.

One can also speak about the principal symbol of an element of $\mathrm{U}(\mathrm{sl}_2)$.
Indeed, $\mathrm{U}(\mathrm{sl}_2)$ is canonically identified with the symmetric
algebra $S(\mathrm{sl}_2)$ as $\mathrm{sl}_2$-modules (see, e.g., \cite[p.82]{dix}).
Using the realization (\ref{gensl}), the algebra $S(\mathrm{sl}_2)$ can be
projected to $\mathrm{Pol}(T^*{{\mathbb R}})$.
Therefore, one can define in a natural way the principal symbol on $S(\mathrm{sl}_2)$.

Our goal is to construct an $\mathrm{sl}_2$-equivariant linear map
$T_{\lambda}:\mathrm{U}(\mathrm{sl}_2) \to {\cal D}_{\lambda}$ which preserves
the principal symbol, i.e., such that the following diagram commutes:
\[
\begin{CD}
\mathrm{U}(\mathrm{sl}_2) @>T_{\lambda} >> {\cal D}_{\lambda} \strut\\
@V{\sigma}VV @VV{\sigma}V \strut\\
\mathrm{Pol}(T^*{{\mathbb R}}) @>id >>
\mathrm{Pol}(T^*{{\mathbb R}})\strut
\end{CD}
\]

{\bf 2.4 Projectively equivariant symbol.}
Viewed as a $\mathrm{Vect}({\mathbb R})$-module, the space of symbols corresponding to
${\cal D}_{\lambda}$ has the form:
\begin{equation}
\label{symb}
\mathrm{Pol}(T^*{{\mathbb R}})
\cong
{\cal F}_0\oplus{\cal F}_1\oplus\cdots\oplus
{\cal F}_n\oplus\cdots.
\end{equation}
The space of polynomials  of degree $\leq n$ is a submodule
of $\mathrm{Pol}(T^*{{\mathbb R}})$ which we denote $\mathrm{Pol}_n(T^*{{\mathbb R}})$.

The following result of \cite{ga}  allows one to identify, for arbitrary $\lambda$,
${\cal D}^n_{\lambda}$ with
$\mathrm{Pol}_n(T^*{{\mathbb R}})$ as $\mathrm{sl}_2$-modules:

(i) There exists a unique  $\mathrm{sl}(2,{\mathbb R})$-isomorphism
$\sigma_\lambda:{\cal D}^n_{\lambda} \to \mathrm{Pol}_n(T^*{{\mathbb R}})$
preserving the principal symbol.

(ii) $\sigma_\lambda$ associates to each differential operator $A$ the polynomial
$\sigma_{\lambda}(A)=\sum\limits_{p=0}^{n} \bar a_p(x) \xi^p $, defined by:
\begin{equation}
\bar a_p(x) = \sum_{j=p}^n\alpha_p^{j}a_j^{(j-p)},
\label{apbar}
\end{equation}
where the constants $\alpha_p^{j}$ are given by:
\[
\alpha_p^{j} =
\frac{{\binom{j}{p}{\binom{2\lambda-p}{j-p}}}}{{\binom{j+p+1}{2p+1}}}
\]
(the binomial coefficient $\binom{\lambda}{j} = \lambda
(\lambda-1) \cdots (\lambda-j+1)/j!$ is a polynomial in
$\lambda$).

The isomorphism $\sigma_\lambda$ is called the \textit{projectively equivariant symbol map}.
Its explicit formula was first found
in~\cite{cmz,lo} in the general case of pseudo-differential operators on
a~one-dimensional manifold
(see also~\cite{lo} for the multi-dimensional case).

\section{Main result}

In this section, we give the main result of this paper.
We adopt the following notations:
\[
{[L_{X_1}^{\lambda}L_{X_2}^{\lambda}\cdots L_{X_n}^{\lambda}]}_+
:=
\sum_{\tau \in S_n} L^{\lambda}_{X_{\tau(1)}}\circ L^{\lambda}_{X_{\tau(2)}}
\circ  \cdots \circ  L^{\lambda}_{X_{\tau(n)}}
\]
for a symmetric $n$-linear map from $\mathrm{Vect}({\mathbb R})$ to $\cal D$
and
\[
{(X_1  X_2 \cdots  X_n)}_{+}
:=
\sum_{\tau \in S_n}
X_{\tau(1)} X_{\tau(2)}  \cdots  X_{\tau(n)}
\]
for a symmetric $n$-linear map from $\mathrm{sl}_2$ to $\mathrm{U}(\mathrm{sl}_2)$,
where $S_n$ is the group of permutations of $n$ elements and
$X_i \in \mathrm{sl}_2$.

\resetfootnoterule

\begin{theorem}
\label{main}
(i) For arbitrary $\lambda \in {{\mathbb C}}$, there exists a unique
$\mathrm{sl}_2$-equivariant linear map preserving the principal symbol:
\[
T_{\lambda}:\mathrm{U}(\mathrm{sl}_2) \to {\cal D}_{\lambda}
\]
defined by
\begin{equation}
T_{\lambda}({(X_1  X_2 \cdots  X_n)}_{+})
=
{[L_{X_1}^{\lambda}L_{X_2}^{\lambda}\cdots L_{X_n}^{\lambda}]}_+  ,
\label{tlform}
\end{equation}
where $X_i \in \{e_0,e_1,e_2 \}$, $L_{X_i}^{\lambda}$
given by (\ref{lieder}) and  $n = 1,2,\dots$.

(ii) The operator $T_{\lambda}$ is given in term of
the $\mathrm{sl}_2$-equivariant symbol (\ref{apbar}) by:
\begin{equation}
\sigma_\lambda([L_{X_1}^{\lambda}L_{X_2}^{\lambda}\cdots L_{X_n}^{\lambda}]_+)
=
\sum_{\mbox{\scriptsize$\begin{matrix} 0\leq k \leq n \\ k \; even\end{matrix}$}}
P^n_k (\lambda)
{\cal A}_k (X_1,\dots,X_n)
\xi^{n-k}  ,
\label{symbT}
\end{equation}
where
\begin{gather}
{\cal A}_k (X_1,\dots,X_n) \nonumber\\
\qquad {}= \sum_{2p+m=k} {\textstyle \binom{k/2}{p}} {(-2)}^p
{(X_1'' \dots X_p'' X_{p+1}' \cdots X_{p+m}'{X}_{p+m+1} \cdots
X_n)}_+ \label{opsy}
\end{gather}
and
\begin{equation}
P^n_k (\lambda)
=
\sum^n_{p=0}
\sum^n_{l=n-k}
 (l-n+k)! \;
\frac{{{\binom{l}{{n-k}}^2}{\binom{2\lambda-n+k}{l-n+k}}}}{{\binom{n-k+l+1}{2n-2k+1}}}
{\textstyle \binom{n}{p}} \mbox{\scriptsize$\left\{\begin{matrix} p \\
l\end{matrix}\right\}$} {\lambda}^{n-p}  , \label{Polsy}
\end{equation}
where $\mbox{\scriptsize$\left\{\begin{matrix} p \\
l\end{matrix}\right\}$}$
 is the Stirling number of the second
 kind\footnote{We refer to \cite{gkp} as a nice elementary introduction to
the combinatorics of the Stirling numbers.}.
\end{theorem}

It is worth noticing that the linear map $T_{\lambda}$ does not depend on the choice
of the PBW-base in $\mathrm{U}(\mathrm{sl}_2)$.

\section{Proof of Theorem \ref{main}}

By construction, the linear map $T_{\lambda}$ is $\mathrm{sl}_2$-equivariant.

{\bf 4.1 $\boldsymbol{\mathrm{sl}_2}$-invariant symmetric differential operators.}
To prove part (ii) of Theorem~\ref{main}
one needs the following
\begin{proposition}
\label{resPro} For arbitrary $\mu \in {{\mathbb C}}$ and $n
=1,2,\dots$, there exists at most one, up to proportionality,
$\mathrm{sl}_2$-equivariant symmetric operator ${\otimes}^n
\mathrm{sl}_2 \to {{\cal F}}_\mu $ which is differential with
respect to the vector fields $X_i \in \mathrm{sl}_2$. This
operator exists if and only if $\mu = k-n$, where~$k$ is an even
positive integer. It is denoted: $ {\cal A}_k  :  {\otimes}^n
\mathrm{sl}_2 \to {{\cal F}}_{k-n}$, and defined by the
expression~(\ref{opsy}).
\end{proposition}

\begin{proof} Each $k$-th order differential operator
${\cal A}:{\otimes}^n \mathrm{sl}_2 \to {{\cal F}}_\mu$ is of the form:
\[
{\cal A} (X_1,\dots,X_n)
= \sum_{2p+m=k}
\beta_p(x)
{(X_1'' \cdots X_p'' X_{p+1}' \cdots X_{p+m}'{X}_{p+m+1} \cdots X_n)}_+,
\]
where $\beta_p(x)$ are some functions.

The condition of $\mathrm{sl}_2$-equivariance for $\cal A$ reads as follows:
\[
X [{\cal A} (X_1,\dots,X_n)]'  +  \mu X' {\cal A} (X_1,\dots,X_n)
=
\sum^n_{i=1}
{\cal A} (X_1,\dots,L_{X}^{-1}(X_i),\dots,X_n)  ,
\]
where $X \in \mathrm{sl}_2$.

Substitute $X=\partial$ to check that the coefficients $\beta_p(x)$
do not depend on $x$. Substitute $X= x\partial$ to obtain the condition $\mu = k-n$.
At last, substitute $X= x^2\partial$ and put $\beta_0 =1$
to obtain, for even $k$, the coefficients from (\ref{opsy}).
If $k$ is odd, one obtains $\beta_p = 0$ for all $p$.

Proposition \ref{resPro} is proven.\end{proof}

The general form (\ref{symbT}) is a consequence of Proposition \ref{resPro}
and decomposition (\ref{symb}).

{\bf 4.2 Polynomials $\boldsymbol{P^n_k(\lambda)}$.}
To compute the polynomials $P^n_k$, put $X_1=\cdots=X_n=x\partial$.
One readily gets, from (\ref{symbT}),
\begin{equation}
\sigma_\lambda ( T_\lambda (X_1,\dots,X_n))|_{x=1}= n!
\sum_{\mbox{\scriptsize $\begin{matrix} 0\leq k \leq n \\ k \; even\end{matrix}$}}
P^n_k (\lambda) \; \xi^{n-k}.
\label{eq1}
\end{equation}
Furthermore, using the well-known expression ${(x\partial)}^n =
\sum\limits_{l=0}^n \mbox{\scriptsize$\left\{\begin{matrix} n\\
l\end{matrix}\right\}$} x^l {\partial}^l$, one has:
\begin{gather*}
T_\lambda (X_1,\dots,X_n)  =  n!  \; (x\partial + \lambda)^n \\
\phantom{T_\lambda (X_1,\dots,X_n) } {} =  n! \sum_{p=0}^n
{\textstyle \binom{n}{p}} {(x\partial)}^n  \lambda^{n-p}
 =  n!  \sum_{p=0}^n \sum^n_{l=0} {\textstyle\binom{n}{p}}
 \mbox{\scriptsize$\left\{\begin{matrix} n\\
 l\end{matrix}\right\}$}
 x^l {\partial}^l \lambda^{n-p} .
\end{gather*}
A straightforward computation
gives the projectively equivariant symbol (\ref{apbar}) of this differential operator:
\begin{gather*}
\sigma_\lambda ( T_\lambda (X_1,\dots,X_n))|_{x=1}\\
\qquad{}=
n! \sum_{\mbox{\scriptsize$\begin{matrix} 0\leq k \leq n \\ k \; even\end{matrix}$}}
\sum^n_{p=0}
\sum^n_{l=n-k}
 (l-n+k)!
\frac{{{\binom{l}{{n-k}}^2}{\binom{2\lambda-n+k}{l-n+k}}}}{{\binom{n-k+l+1}{2n-2k+1}}}
{\textstyle\binom{n}{p}} \mbox{\scriptsize$\left\{\begin{matrix} p\\
l\end{matrix}\right\}$} {\lambda}^{n-p} \xi^{n-k} .
\end{gather*}
Compare with the equality (\ref{eq1}) to obtain the formulae from (\ref{Polsy}).

Theorem \ref{main} (ii) is proven.

{\bf 4.3 Uniqueness.}
Let $T$ be an $\mathrm{sl}_2$-equivariant linear map
$\mathrm{U}(\mathrm{sl}_2) \to {\cal D}_{\lambda}$
for a certain $\lambda \in {{\mathbb C}}$.
In view of the decomposition (\ref{symb}), it follows from Proposition \ref{resPro} that
$\sigma_\lambda \circ T |_{{\cal F}_k}=c_k(\lambda) {\cal A}_k$,
where $c_k(\lambda)$ is a constant depending on $\lambda$.
Recall that $\mathrm{Pol}_n(T^*{{\mathbb R}})$ is a {\it rigid}
$\mathrm{sl}_2$-module, i.e., every $\mathrm{sl}_2$-equivariant
linear map on $\mathrm{Pol}_n(T^*{{\mathbb R}})$ is proportional to the identity
(see, e.g., \cite{lo}).
Assuming, now, that $T$ preserves the principal symbol, the rigidity
of $\mathrm{Pol}_n(T^*{{\mathbb R}})$
fixes the constants $c_k(\lambda)$ in a unique way. Hence the uniqueness of $T_\lambda$.

Theorem \ref{main} is proven.

\section{The embedding $\boldsymbol{{\rm gl(\lambda)} \to {\cal D}_{\lambda}}$}

A corollary of the uniqueness of the operator $T_\lambda$ and results
of \cite{bb1,bb2,fe,sh}
is that the embedding ${\rm gl(\lambda)} \to {\cal D}_{\lambda}$
constructed in \cite{fe} coincides with $T_\lambda$.

More precisely, according to results of \cite{bb1,bb2,sh},
there exists a homomorphism of Lie algebras
$p_\lambda : \mathrm{U}(\mathrm{sl}_2) \to {\cal D}_{\lambda}$
preserving the principal symbol.
The homomorphism $p_\lambda$ is, in particular, $\mathrm{sl}_2$-equivariant. By uniqueness
of $T_\lambda$, one has $T_\lambda=p_\lambda$. It is also proven
that the  kernel of $p_\lambda$ is a two-sided
ideal of $\mathrm{U}(\mathrm{sl}_2)$ generated
by $\Delta - \lambda(\lambda-1)$ (see \cite{bb1,bb2}).
Taking the quotient, one then has an
embedding ${\tilde T_\lambda}:{\rm gl(\lambda)} \to {\cal D}_{\lambda}$.
Since the embedding from~\cite{fe} preserves the principal symbol,
it is equal to ${\tilde T_\lambda}$. Finally, it is obvious that
the image of $T_\lambda$ is the subalgebra
${\cal D}^{\rm pol}_{\lambda} \subset {\cal D}_{\lambda} $
of differential operators with polynomial coefficients.
Therefore,
${\tilde T_\lambda}:{\rm gl(\lambda)} \to {\cal D}^{\rm pol}_{\lambda}$
is a Lie algebras isomorphism.

\section{Examples}

As an illustration of Theorem \ref{main}, let us give the expressions of the general
formulae~(\ref{tlform}) and (\ref{symbT}) for the order $n=1,2,3,4,5$.
Let $X_1,X_2,X_3,X_4$ and $X_5$ be arbitrary vector fields in $\mathrm{sl}_2$.

1) The $\mathrm{sl}_2$-equivariant symbol,
defined by (\ref{apbar}), of a first order operator of a Lie derivative
$L_{X_1}^{\lambda}$ is
\[
\sigma_\lambda (L_{X_1}^{\lambda})  =  X_1 (x) \xi.
\]

2) The ``anti-commutator'' ${[L_{X_1}^{\lambda}L_{X_2}^{\lambda}]}_+$
has the following projectively equivariant symbol:
\[
\sigma_\lambda ({[L_{X_1}^{\lambda}L_{X_2}^{\lambda}]}_+)=
{({X_1}{X_2})}_+ \xi^2
+\frac{1}{3}\lambda(\lambda-1)({({X_1'}{X_2'})}_+ -2 {(X_1'' X_2)}_+)
\]
which also following from (\ref{apbar}).

3) The projectively equivariant symbol of a third order expression
${[L_{X_1}^{\lambda}L_{X_2}^{\lambda}L_{X_3}^{\lambda}]}_+$
can be also easily calculated from (\ref{apbar}). The result is:
\begin{gather*}
\sigma_\lambda({[L_{X_1}^{\lambda}L_{X_2}^{\lambda}L_{X_3}^{\lambda}]}_+)
 ={({X_1}{X_2}{X_3})}_+ \xi^3 \\
\qquad{}+\frac{1}{5}(3\lambda^2-3\lambda-1)
({(X_1' X_2' X_3)}_+ -2{(X_1'' X_2 X_3)}_+) \xi .
\end{gather*}

4) Direct calculation from (\ref{apbar}) gives the
projectively equivariant symbol of a fourth order expression
${[L_{X_1}^{\lambda}L_{X_2}^{\lambda}L_{X_3}^{\lambda}L_{X_4}^{\lambda}]}_+$, that is:
\begin{gather*}
\sigma_\lambda({[L_{X_1}^{\lambda}L_{X_2}^{\lambda}L_{X_3}^{\lambda}L_{X_4}^{\lambda}]}_+)
={({X_1}{X_2}{X_3}{X_4})}_+  \xi^4 \\
\qquad{}+\frac{1}{7}(6\lambda^2-6\lambda-5)
({(X_1' X_2' X_3 X_4)}_+ -2{(X_1'' X_2 X_3 X_4)}_+)  \xi^2 \\
\qquad {}+\frac{1}{15}\lambda(\lambda-1)(3\lambda^2-3\lambda-1)
({(X_1' X_2' X_3' X_4')}_+ -4{(X_1'' X_2' X_3' X_4)}_+ \\
\qquad{}+4{(X_1'' X_2'' X_3 X_4)}_+)  .
\end{gather*}

5) In the same manner, one can easily check that the $\mathrm{sl}_2$-equivariant symbol of
a fifth order expression
${[L_{X_1}^{\lambda}
L_{X_2}^{\lambda}L_{X_3}^{\lambda}L_{X_4}^{\lambda}L_{X_4}^{\lambda}]}_+$ is:
\begin{gather*}
\sigma_\lambda({[L_{X_1}^{\lambda}L_{X_2}^{\lambda}L_{X_3}^{\lambda}
L_{X_4}^{\lambda}L_{X_5}^{\lambda}]}_+)
={({X_1}{X_2}{X_3}{X_4}{X_5})}_+  \xi^5 \\
\qquad{}+\frac{5}{9}(2\lambda^2-2\lambda-3)
({(X_1' X_2' X_3 X_4 X_5)}_+ -2{(X_1'' X_2 X_3 X_4 X_5)}_+)  \xi^3\\
\qquad {}+\frac{1}{7}(3\lambda^4-6\lambda^3+3\lambda+1)
({(X_1' X_2' X_3' X_4' X_5)}_+ -4{(X_1'' X_2' X_3' X_4 X_5)}_+\\
\qquad{}+4{(X_1'' X_2'' X_3 X_4 X_5)}_+)\xi .
\end{gather*}

\subsection*{Acknowledgments}
I would like to thank V~Ovsienko for statement of
the problem.
I am also grateful to  Ch~Duval and A~El~Gradechi
for enlightening  discussions.

\label{Gargoubi-lastpage}

\begin{thebibliography}{99}\small

\bibitem{bb1} Beilinson~A and Bernstein~J, Localisation de $\mathrm{g}$-modules,
{\it C.R. Acad. Sci. Paris Ser. I Math.} {\bf 292} (1981), 15--18.

\bibitem{bb2} Beilinson A and Bernstein J, A Proof of Jantzen Conjectures,
{\it Adv. in Sov. Math.} {\bf 16} (1993), 1--50.

\bibitem{car} Cartan E, Le\c cons sur la th\'eorie des espaces
\`a connexion projective, Gauthier -- Villars, Paris, 1937.

\bibitem{cmz} Cohen P, Manin Yu and Zagier D,
Automorphic Pseudodifferential Operators, in Progr. Nonlinear Diff.
Eq. Appl., Vol.~26, Birkh\"auser, Boston, 1997, 17--47.

\bibitem{do}
Duval C and Ovsienko~V, Space of Second Order Linear Differential
Operators as a Module Over the Lie Algebra of Vector Fields, {\it
Adv. in Math.} {\bf 132}, Nr.~2 (1997), 316--333.

\bibitem{dix}
Dixmier J,
Alg\`ebres enveloppantes, Gauthier -- Villars, Paris, 1974.

\bibitem{fe}
Feigin B~L, The Lie Algebras ${\rm gl}(\lambda)$
and Cohomologies of Lie Algebra of Differential
Operators, {\it Russian Math. Surveys},
{\bf 43}, Nr.~2 (1988), 157--158.

\bibitem{ga} Gargoubi H,
Sur la g\'eom\'etrie de l'espace des op\'erateurs diff\'erentiels
lin\'eaires sur $\bf R$, {\it Bull. Soc. Roy. Sci. Li\`{e}ge} {\bf
69}, Nr.~1 (2000), 21--47.

\bibitem{go} Gargoubi H and Ovsienko V,
Space of Linear Differential Operators on the Real Line as a
Module Over the  Lie Algebra of Vector Fields, {\it Internat.
Mathem. Res. Notices}~Nr.~5 (1996), 235--251.

\bibitem{go1} Gargoubi H and Ovsienko V,
Modules of Differential Operators on the Real Line, {\it Funct.
Anal. Appl.} {\bf 35}, Nr.~1 (2001), 16--22.

\bibitem{gkp} Graham R, Knuth D and Patashnik O,
Concrete Mathematics,
Addison-Wesley, 1989.

\bibitem{gl}
Grozman P and Leites D~A,
Lie Superalgebras of Supermatrices of Complex Size. Their Generalizations and
Related Integrable Systems, in  Proc. Internatnl.
Symp. Complex Analysis and Related Topics,
Editors: E. Ramirez de Arellano, et. al., Mexico, 1996, Birkh\"auser Verlag, 1999, 73--105.


\bibitem{km}
Khesin B and Malikov F,
Universal Drinfeld--Sokolov Reduction and the Lie Algebras of Matrices of Complex Size,
{\it Comm. Math. Phys.} {\bf 175}, Nr.~1 (1996), 113--134.


\bibitem{lmt}
Lecomte P~B~A, Mathonet~P and Tousset~E, Comparison of Some
Modules of the Lie Algebra of Vector Fields, {\it Indag. Math.,
N.S.} {\bf 7}, Nr.~4 (1996), 461--471.


\bibitem{lo}
Lecomte P~B~A and Ovsienko~V,
Projectively Invariant Symbol Calculus, {\it Lett. Math. Phys.} {\bf 49}, Nr.~3 (1999),
173--196.

\bibitem{ls}
Leites D~A and Sergeev A~N,
Orthogonal Polynomials of  a Discrete Variable and Lie Algebras of Complex-Size Matrices,
{\it Theor. Math. Phys.} {\bf 123}, Nr.~2 (2000), 582--608.

\bibitem{sh} Shoikhet B, Certain Topics on the Representation Theory of the Lie Algebra
${\rm gl}(\lambda)$. Complex Analysis and Representation Theory.~1,
{\it J. Math. Sci. (New York)} {\bf 92}, Nr.~2
(1998), 3764--3806.

\bibitem{wil} Wilczynski E~J, Projective Differential Geometry
of Curves and Ruled Surfaces, Leipzig -- Teubner, 1906.
\end{thebibliography}
\end{document}